\documentclass[12pt,reqno]{amsart}

\usepackage{amsmath, epsfig, cite}
\usepackage{amssymb}
\usepackage{amsfonts}
\usepackage{latexsym}
\usepackage{graphicx}
\usepackage{algorithmic,algorithm}
\usepackage{amsthm,array}

\newtheorem{thm}{Theorem}[section]

\newtheorem{lem}[thm]{Lemma}

\newtheorem{exam}[thm]{Example}
\newtheorem{rem}[thm]{Remark}

\def\pf{\noindent {\it Proof.} }

\numberwithin{equation}{section}

\makeatletter \@addtoreset{equation}{section} \makeatother

\begin{document}
\allowdisplaybreaks
\title
{Arithmetic properties of the $\ell$-regular partitions}
\author[Suping Cui]{Suping Cui}
\address[S. P. Cui]{Center for Combinatorics, LPMC-TJKLC, Nankai
University, Tianjin 300071, P. R. China} \email{jiayoucui@163.com}

\author[N. S. S. Gu]{Nancy Shanshan Gu}
\address[N. S. S. Gu]{Center for Combinatorics, LPMC-TJKLC, Nankai
University, Tianjin 300071, P. R. China} \email{gu@nankai.edu.cn}

\date{\today}
\begin{abstract} For a given prime $p$, we study the properties of the $p$-dissection identities
of Ramanujan's theta functions $\psi(q)$ and $f(-q)$, respectively. Then as applications, we find many infinite family of congruences modulo $2$
for some $\ell$-regular partition functions, especially, for $\ell=2,4,5,8,13,16$.
Moreover, based on the classical congruences for $p(n)$ given by Ramanujan,
we obtain many more congruences for some $\ell$-regular partition functions.

\end{abstract}

\maketitle

\section{Introduction}

A partition of a positive integer $n$ is a nonincreasing sequence of positive
integers whose sum is $n$. Let $p(n)$ denote the number of partitions
of $n$. If $\ell$ is a positive integer, then a partition is called
$\ell$-regular partition if there is no part divisible by $\ell$.
Let $b_\ell(n)$ denote the number of $\ell$-regular partitions of
$n$. The generating function of $b_\ell(n)$ is stated as
follows.
\begin{equation*}
\sum_{n=0}^{\infty}b_\ell(n)q^n=\frac{(q^\ell;q^\ell)_{\infty}}{(q;q)_{\infty}}.
\end{equation*}

The divisibility and distribution of $b_\ell(n)$ modulo $m $ were
widely studied in the literature, see
\cite{Ahlgren-Lovejoy-2001,Calkin-Drake-James-2008,Chen-2011,Dandurand-Penniston-2009,Furcy-Penniston-2012,Gordon-Ono-1997,Hirschhorn-Sellers-2010,
Andrews-Hirschhorn-Sellers-2010,Lovejoy-2001,Lovejoy-2003,Lovejoy-Penniston-2001,Ono-Penniston-2000,Penniston-2002,Penniston-2008,Webb-2011}.
Recently, the arithmetic of the $\ell$-regular partitions has received
a great deal of attentions. For example, in \cite[Theorem 3.5]{Andrews-Hirschhorn-Sellers-2010}, Andrews et al. gave
some infinite family congruences modulo $2$ for the $4$-regular
partition function. They used $ped(n)$ to denote the number of the $4$-regular partitions of $n$ which were called partitions with even parts distinct.
For $\alpha \geq 0$ and $n \geq 0$,
\begin{align}
ped(3^{2\alpha+1}n+\frac{17\cdot 3^{2\alpha}-1}{8}) & \equiv 0
\quad (\text{mod}\ 2),\label{b4-inf-23}\\
ped(3^{2\alpha+2}n+\frac{11 \cdot 3^{2\alpha+1}-1}{8}) & \equiv 0
\quad (\text{mod}\ 2),\label{b4-inf-24}\\
ped(3^{2\alpha+2}n+\frac{19 \cdot 3^{2\alpha+1}-1}{8}) & \equiv 0
\quad (\text{mod}\ 2).\label{b4-inf-25}
\end{align}

In \cite{Chen-2011}, Chen used the theory of Hecke eigenforms to get more generalized congruences modulo $2$ for $b_4(n)$.

In \cite{Calkin-Drake-James-2008}, Calkin et al. studied the
$2$-divisibility of $b_5(n)$, and gave the following congruences.
\begin{thm}\label{Calkin-b5}\cite[Theorem 3]{Calkin-Drake-James-2008} For $n \geq 0$,
\begin{align*}
b_5(20n+5) \equiv 0 \quad (\text{mod}\ 2),\\
b_5(20n+13) \equiv 0 \quad (\text{mod}\ 2).
\end{align*}
\end{thm}
Later, Hirschhorn and Sellers got the congruence for $b_5(n)$ in \cite{Hirschhorn-Sellers-2010} .
\begin{thm}\cite[Theorem 2.5]{Hirschhorn-Sellers-2010} \label{Hir-b5}Suppose that $p$ is any prime greater than $3$ such that $-10$ is
a quadratic nonresidue modulo $p$, $u$ is the reciprocal of $24$ modulo $p^2$, and $r \not\equiv 0\ (\text{mod}\ p)$. Then, for all $m$,
\begin{equation*}
b_5(4p^2m+4u(pr-7)+1)\equiv0\quad (\text{mod}\ 2).
\end{equation*}
\end{thm}

In \cite{Calkin-Drake-James-2008}, Calkin et al. posed the following conjecture for $b_{13}(n)$ which was proved by Webb in
\cite{Webb-2011}.

\begin{thm}\cite{Webb-2011} For $n\geq 0$ and $\ell\geq 2$, there holds
\begin{align}\label{133333}
b_{13}(3^{\ell}n+\frac{5\cdot 3^{\ell-1}-1}{2}) \equiv 0 \quad
(\text{mod}\ 3).
\end{align}
\end{thm}

By using a modification of the method which was utilized in \cite{Webb-2011}, Furcy and Penniston
obtained the following congruences modulo $3$ for some $\ell$-regular partition functions in \cite{Furcy-Penniston-2012}.
\begin{thm}\cite[Theorem 1]{Furcy-Penniston-2012}\label{theom1.4}
For $\alpha\geq 0$ and $n\geq 0$, we have
\begin{align*}
&b_{7}(3^{2\alpha+2}n+\frac{11\cdot 3^{2\alpha+1}-1}{4})\equiv b_{7}(3^{2\alpha+3}n+\frac{5\cdot 3^{2\alpha+2}-1}{4}) \equiv 0 \quad (\text{mod}\ 3),\\
&b_{19}(3^{2\alpha+4}n+\frac{5\cdot 3^{2\alpha+3}-3}{4})\equiv b_{19}(3^{2\alpha+5}n+\frac{11\cdot 3^{2\alpha+4}-3}{4}) \equiv 0 \quad (\text{mod}\ 3),\\
&b_{25}(3^{2\alpha+3}n+2\cdot 3^{2\alpha+2}-1) \equiv 0 \quad (\text{mod}\ 3),\\
&b_{34}(3^{4\alpha+3}n+\frac{19\cdot 3^{4\alpha+2}-11}{8})\equiv b_{34}(3^{4\alpha+5}n+\frac{11\cdot 3^{4\alpha+4}-11}{8}) \equiv 0 \quad (\text{mod}\ 3),\\
&b_{37}(3^{3\alpha+3}n+\frac{3^{3\alpha+2}-3}{2}) \equiv 0 \quad (\text{mod}\ 3),\\
&b_{43}(3^{2\alpha+4}n+\frac{5 \cdot 3^{2\alpha+3}-7}{4})\equiv b_{43}(3^{2\alpha+5}n+\frac{11\cdot 3^{2\alpha+4}-7}{4}) \equiv 0 \quad (\text{mod}\ 3),\\
&b_{49}(3^{3\alpha+3}n+2\cdot 3^{3\alpha+2}-2) \equiv 0 \quad (\text{mod}\ 3).
\end{align*}
\end{thm}

\begin{thm}\cite[Theorem 2]{Furcy-Penniston-2012}\label{fp-t2} For every $n \geq 0$,
\begin{equation*}
b_{10}(9n+3)\equiv b_{22}(27n+16)\equiv b_{28}(27n+9)\equiv 0 \quad (\text{mod}\ 3).
\end{equation*}
\end{thm}

During the study of the congruences for the $\ell$-regular partition functions, we observe that some
generating functions of the $\ell$-regular partitions congruent to the functions related to Ramanujan's theta functions $\psi(q)$ and $f(-q)$ modulo $2$.
Then, for a given prime $p$, by finding the properties of the $p$-dissection identities of $\psi(q)$ and $f(-q)$, we obtain
many infinite family congruences modulo $2$ for $b_{\ell}(n)$, such as $l=2,4,5,8,13,16$.
Especially, we generalize the congruences modulo 2 for $b_4(n)$ given by Chen in \cite{Chen-2011} which include \eqref{b4-inf-23}-\eqref{b4-inf-25}
given by Andrews et al. in \cite{Andrews-Hirschhorn-Sellers-2010} as the special cases, and we also generalize the
congruences modulo $2$ for $b_5(n)$ in Theorem \ref{Calkin-b5} given by Calkin et al. in \cite{Calkin-Drake-James-2008} and in Theorem \ref{Hir-b5} given by Hirschhorn and Sellers in \cite{Hirschhorn-Sellers-2010}.
Meanwhile, by observing a relation between the generating function of $b_4(n)$ and that of $b_{13}(n)$, we get some new congruences for $b_{13}(n)$.
Moreover, by combining the classical congruences for $p(n)$ given by Ramanujan with some congruences obtained in this paper and some
known congruences in Theorem \ref{theom1.4} and Theorem \ref{fp-t2}, we get many more congruences for some $\ell$-regular partition functions.

In section 2, we study the $p$-dissection identities of $\psi(q)$ and $f(-q)$,
respectively. In section 3, we mainly focus on finding the applications
of this two identities on the arithmetic properties of the $\ell$-regular partitions modulo $2$ for $\ell=2,4,5,8,13,16$.
In section 4, based on Ramanujan's famous congruences for $p(n)$, we obtain some more congruences for the $\ell$-regular partition functions.

As usual, we follow the notation and terminology in
\cite{Gasper-Rahman-2004}. For $|q|<1$, the $q$-shifted factorial is
defined by
\begin{equation*}
(a;q)_\infty= \prod_{k=0}^{\infty}(1-aq^k) \quad\text{and}\quad
(a;q)_n = \frac{(a;q)_\infty}{(aq^n;q)_\infty}, \text{ for } n\in
\mathbb{C}.
\end{equation*}

In the following, we list some definitions and identities which are
frequently used in this paper.

The Legendre symbol is a function of $a$ and $p$ defined as follows:
$$\left(\frac{a}{p}\right)=\left\{\begin{array}{ll}1,&\text{if } a \text{ is a quadratic residue modulo } p \text{ and } a \not\equiv 0\ (\text{mod}\ p),\\
-1,&\text{if } a \text{ is a quadratic non-residue modulo } p,\\
0,&\text{if } a\equiv 0\ (\text{mod}\ p). \end{array}\right.$$

Jacobi's triple product identity \cite[Theorem 1.3.3]{Berndt-2004}:
for $z \neq 0$ and $|q|<1$,
\begin{equation}\label{Jacobi}
\sum_{n=-\infty}^{\infty}z^nq^{n^2}=(-zq,-q/z,q^2;q^2)_{\infty}.
\end{equation}

Euler's pentagonal number theorem \cite[Corollary
1.3.5]{Berndt-2004}:
\begin{equation}\label{Euler}
\sum_{n={-\infty}}^{\infty}(-1)^nq^{\frac{n(3n+1)}{2}}=(q;q)_{\infty}.
\end{equation}

Jacobi's identity \cite[Theorem 1.3.9]{Berndt-2004}:
\begin{equation}\label{Jacobi-3}
\sum_{n=0}^{\infty}(-1)^n(2n+1)q^{\frac{n(n+1)}{2}}=(q;q)_{\infty}^3.
\end{equation}

Ramanujan's general theta function $f(a,b)$ is defined by
\begin{equation*}
f(a,b):=\sum_{n=-\infty}^{\infty}a^{\frac{n(n+1)}{2}}b^{\frac{n(n-1)}{2}},\qquad
|ab|<1.
\end{equation*}
The special cases of $f(a,b)$ are stated as follows.
\begin{align*}
\psi(q)&:=f(q,q^3)=\sum_{n=0}^{\infty}q^{\frac{n(n+1)}{2}}=\frac{(q^2;q^2)_{\infty}}{(q;q^2)_{\infty}},\\
f(-q)&:= f(-q,-q^2)=\sum_{n=-\infty}^{\infty}(-1)^nq^{\frac{n(3n-1)}{2}}=(q;q)_{\infty}.
\end{align*}

The quintuple product identity \cite[Theorem 1.3.19]{Berndt-2004}:
\begin{equation}\label{quintuple}
\frac{f(-x^2,-\lambda x)f(-\lambda x^3)}{f(-x,-\lambda x^2)}=f(-\lambda^2x^3,-\lambda x^6)+xf(-\lambda,-\lambda^2x^9).
\end{equation}

\section{Preliminaries}

In this section, we give a $p$-dissection identity of $\psi(q)$ for any odd prime $p$,
and a $p$-dissection identity of $f(-q)$ for any prime $p \geq 5$.

\begin{thm}\label{psi-p}
For any odd prime $p$, we have
\begin{equation*}
\psi(q)
=\sum_{k=0}^{\frac{p-3}{2}}q^{\frac{k^{2}+k}{2}}f(q^{\frac{p^2+(2k+1)p}{2}},q^{\frac{p^2-(2k+1)p}{2}})+q^{\frac{p^{2}-1}{8}}\psi(q^{p^2}).
\end{equation*}
Meanwhile, we claim that $(k^2+k)/2$ and $(p^2-1)/8$ are
not in the same residue class modulo $p$ for $0 \leq k \leq (p-3)/2$.

\end{thm}
\pf We have
\begin{equation*}
\psi(q)=\sum_{n=0}^{\infty}q^{\frac{n(n+1)}{2}}=\frac{1}{2}\sum_{n=-\infty}^{\infty}q^{\frac{n(n+1)}{2}}.
\end{equation*}
For any odd prime $p$, we separate out the above series terms in
powers of $q$ according to their residue classes modulo $p$. Then,
we have
\begin{align}
\psi(q) &
=\frac{1}{2}\sum_{k=0}^{p-1}\sum_{n={-\infty}}^{\infty}q^{\frac{(pn+k)(pn+k+1)}{2}} \nonumber\\
&=\frac{1}{2}\sum_{k=0}^{p-1}\sum_{n={-\infty}}^{\infty}q^{\frac{p^{2}n^{2}+(2k+1)pn+k^{2}+k}{2}}\nonumber\\
&=\frac{1}{2}\sum_{k=0}^{\frac{p-3}{2}}\sum_{n={-\infty}}^{\infty}q^{\frac{p^{2}n^{2}+(2k+1)pn+k^{2}+k}{2}}
+\frac{1}{2}q^{\frac{p^2-1}{8}}\sum_{n={-\infty}}^{\infty}q^{p^2\frac{n(n+1)}{2}}\nonumber\\
&\quad+\frac{1}{2}\sum_{k=\frac{p+1}{2}}^{p-1}\sum_{n={-\infty}}^{\infty}q^{\frac{p^{2}n^{2}+(2k+1)pn+k^{2}+k}{2}} \label{kn}\\
&=\sum_{k=0}^{\frac{p-3}{2}}q^{\frac{k^{2}+k}{2}}\sum_{n={-\infty}}^{\infty}
q^{\frac{p^{2}n^{2}+(2k+1)pn}{2}}
+q^{\frac{p^{2}-1}{8}}\psi(q^{p^2})\nonumber\\
&=\sum_{k=0}^{\frac{p-3}{2}}q^{\frac{k^{2}+k}{2}}f(q^{\frac{p^2+(2k+1)p}{2}},q^{\frac{p^2-(2k+1)p}{2}})+q^{\frac{p^{2}-1}{8}}\psi(q^{p^2})
\nonumber \qquad \text{by\ }\eqref{Jacobi}.
\end{align}
Note that the second summation in \eqref{kn} is the case
$k=(p-1)/2$. Then setting $k \rightarrow p-1-k$ and $n
\rightarrow -n-1$ in the third summation in \eqref{kn}, we get the same expression as
the first summation.

For $0 \leq k,m \leq (p-1)/2$ and $k \neq m$, we know that $0<|k-m|\leq (p-1)/2$ and $1<k+m+1<p$. Then we have
\begin{equation}\label{km}
(k-m)(k+m+1) \not\equiv 0 \quad (\text{mod}\ p).
\end{equation}
So
\begin{equation*}
\frac{k^2+k}{2} \not\equiv \frac{m^2+m}{2} \quad (\text{mod}\ p).
\end{equation*}

Therefore, we know that $(k^2+k)/2$ and $(p^2-1)/8$ can not be in the same residue class modulo $p$ for $0 \leq k \leq (p-3)/2$.
\qed

In the next, we dicuss a $p$-dissection identity of $f(-q)$.
\begin{thm}\label{f-p}For any prime $p \geq 5$, we have
\begin{equation*}
f(-q)=\sum_{\tiny \begin{array}{l}k=-\frac{p-1}{2}\\k\neq \frac{\pm p-1}{6}\end{array}}^{\frac{p-1}{2}} (-1)^kq^{\frac{3k^{2}+k}{2}}f(-q^{\frac{3p^2+(6k+1)p}{2}},-q^{\frac{3p^2-(6k+1)p}{2}})+(-1)^{\frac{\pm p-1}{6}}q^{\frac{p^2-1}{24}}f(-q^{p^{2}}),
\end{equation*}
where $\pm$ depends on the condition that $(\pm p-1)/6$ should
be an integer. Meanwhile, we claim that $(3k^2+k)/2$ and
$(p^2-1)/24$ are not in the same residue class modulo $p$ for $-(p-1)/2
\leq k \leq (p-1)/2$ and $k \neq (\pm p-1)/6$.
\end{thm}
\pf For $p \geq 5$, we have
\begin{align*}
f(-q) &=\sum_{n=-\infty}^{\infty}(-1)^nq^{\frac{n(3n+1)}{2}}\\
&=\sum_{k=-\frac{p-1}{2}}^{\frac{p-1}{2}}\sum_{n={-\infty}}^{\infty}(-1)^{pn+k}q^{\frac{(pn+k)(3pn+3k+1)}{2}} \\
&=\sum_{k=-\frac{p-1}{2}}^{\frac{p-1}{2}}(-1)^k q^{\frac{3k^2+k}{2}}\sum_{n={-\infty}}^{\infty}(-1)^n q^{\frac{3p^2 n^2+(6k+1)pn}{2}}\\
&=\sum_{\tiny \begin{array}{l}k=-\frac{p-1}{2}\\k\neq \frac{\pm p-1}{6}\end{array}}^{\frac{p-1}{2}}(-1)^k q^{\frac{3k^2+k}{2}}
\sum_{n={-\infty}}^{\infty}(-1)^n q^{\frac{3p^2 n^2+(6k+1)pn}{2}}\\
&\qquad\qquad\qquad\qquad\qquad\qquad
+(-1)^{\frac{\pm p-1}{6}}q^{\frac{p^2-1}{24}}\sum_{n={-\infty}}^{\infty}(-1)^nq^{p^2\frac{n(3n+1)}{2}}\\
&=\sum_{\tiny \begin{array}{l}k=-\frac{p-1}{2}\\k\neq \frac{\pm p-1}{6}\end{array}}^{\frac{p-1}{2}}
(-1)^kq^{\frac{3k^{2}+k}{2}}f(-q^{\frac{3p^2+(6k+1)p}{2}},-q^{\frac{3p^2-(6k+1)p}{2}})
+(-1)^{\frac{\pm p-1}{6}}q^{\frac{p^2-1}{24}}f(-q^{p^{2}}).
\end{align*}

Here $\pm$ in $(\pm p-1)/6$ depends on the condition that $(\pm p-1)/6$ should
be an integer.

In addition, set $m=(\pm p-1)/6$. Then for any integer $k$ with $-(p-1)/2
\leq k \leq (p-1)/2$, if we have
\begin{equation*}
\frac{3k^2+k}{2} \equiv \frac{3m^2+m}{2}\quad (\text{mod}\ p),
\end{equation*}
then
\begin{equation*}
(k-m)(3k+3m+1)\equiv 0 \quad (\text{mod}\ p).
\end{equation*}
It is obvious that $k-m \not\equiv 0\ (\text{mod}\ p)$, so
\begin{align*}
3k+3m+1 &\equiv 0 \quad (\text{mod}\ p),\\
6k\pm p +1 &\equiv 0 \quad (\text{mod}\ p),\\
6k+1 &\equiv 0 \quad (\text{mod}\ p).
\end{align*}
Since $-3p+4 \leq 6k+1 \leq 3p-2$ and $6k+1$ is odd, we have
$$6k+1=\pm p.$$
So the only solution is $k=m$. Therefore, for $-(p-1)/2
\leq k \leq (p-1)/2$ and $k\neq m$, we have
\begin{equation*}
\frac{3k^2+k}{2} \not\equiv \frac{3m^2+m}{2}\quad (\text{mod}\ p).
\end{equation*}
\qed

In the following, we consider a special case of Theorem \ref{f-p} by setting $p=5$. Using the quintuple product identity \eqref{quintuple},
we can easily obtain the following identity given by Ramanujan in \cite[p. 212]{Ramanujan-2000}
\begin{equation}\label{Ramanujan-qq}
(q;q)_{\infty}=\frac{(q^{10},q^{15},q^{25};q^{25})_{\infty}}{(q^5,q^{20};q^{25})_{\infty}}-q(q^{25};q^{25})_{\infty}
-q^2\frac{(q^5,q^{20},q^{25};q^{25})_{\infty}}{(q^{10},q^{15};q^{25})_{\infty}}.
\end{equation}
Recently, Hirschhorn gave a simple proof of the above identity by using Jacobi's triple product identity \eqref{Jacobi} in \cite{Hirschhorn-2011}.
We know that this identity plays an important role in the proof of Ramanujan's ``most beautiful identity"
\begin{equation*}
\sum_{n=0}^{\infty}p(5n+4)q^n=5\frac{(q^5;q^5)_{\infty}^5}{(q;q)_{\infty}^6}.
\end{equation*}

\section{Arithmetic of the $\ell$-regular parititions related to $\psi(q)$ and $f(-q)$}

In this section, based on Theorem \ref{psi-p} and Theorem \ref{f-p},
we study some infinite family of congruences for some $\ell$-regular
partition functions for $\ell=2,4,5,8,13,16$.

\subsection{$2$-regular partitions}
In the literature, the $2$-regular partitions are usually called distinct partitions. Combining the following fact
\begin{equation}\label{b2-f(-q)}
\sum_{n=0}^\infty b_2(n)q^n =(-q;q)_{\infty} \equiv f(-q) \quad (\text{mod}\ 2)
\end{equation}
with the $p$-dissection identity of $f(-q)$ in Theorem \ref{f-p}, we find some Ramanujan-type
congruences modulo $2$ for $b_2(n)$.

\begin{lem}\label{b2-f} For any prime $p \geq 5$, $\alpha \geq 0$, and $n \geq 0$, we have
\begin{equation*}
\sum_{n=0}^{\infty}b_{2}(p^{2\alpha}n+\frac{p^{2\alpha}-1}{24})q^{n}\equiv f(-q) \quad (\text{mod}\ 2).
\end{equation*}
\end{lem}
\pf
We prove the lemma by induction on $\alpha$.

Note that \eqref{b2-f(-q)} is the case for $\alpha=0$.

Suppose that the lemma holds for $\alpha$. Then we have
\begin{equation*}
\sum_{n=0}^{\infty}b_{2}(p^{2\alpha}n+\frac{p^{2\alpha}-1}{24})q^{n}\equiv f(-q) \quad (\text{mod}\ 2).
\end{equation*}
According to Theorem \ref{f-p}, we have
\begin{equation}\label{b2-fp}
\sum_{n=0}^{\infty}b_{2}(p^{2\alpha}(pn+\frac{p^2-1}{24})+\frac{p^{2\alpha}-1}{24})q^{n} \equiv f(-q^p) \quad (\text{mod}\ 2).
\end{equation}
Then
\begin{align*}
\sum_{n=0}^{\infty}b_{2}(p^{2\alpha}(p^2n+\frac{p^2-1}{24})+\frac{p^{2\alpha}-1}{24})q^{n}&=\sum_{n=0}^{\infty}b_{2}(p^{2\alpha+2}n+\frac{p^{2\alpha+2}-1}{24})q^{n}
\\ &\equiv f(-q) \quad (\text{mod}\ 2).
\end{align*}
Therefore, the lemma holds for $\alpha+1$.\qed

According to Lemma \ref{b2-f}, we get the following congruences for $b_2(n)$.
\begin{thm}\label{theo1}For any prime $p \geq 5$, $\alpha \geq 1$, and $n \geq 0$, we have
\begin{equation*}
b_{2}(p^{2\alpha}n+\frac{(24i+p)p^{2\alpha-1}-1}{24}) \equiv 0 \quad (\text{mod}\ 2),
\end{equation*}
where $i=1,2,\ldots,p-1$.
\end{thm}
\pf According to \eqref{b2-fp}, we have
\begin{equation*}
\sum_{n=0}^{\infty}b_{2}(p^{2\alpha+1}n+\frac{p^{2\alpha+2}-1}{24})q^{n} \equiv f(-q^p) \quad (\text{mod}\ 2).
\end{equation*}
Therefore, for $i=1,2,\ldots,p-1$, we have
\begin{equation*}
b_{2}(p^{2\alpha+1}(pn+i)+\frac{p^{2\alpha+2}-1}{24}) \equiv 0 \quad (\text{mod}\ 2).
\end{equation*}
That is
\begin{equation*}
b_{2}(p^{2\alpha+2}n+\frac{(24i+p)p^{2\alpha+1}-1}{24}) \equiv 0 \quad (\text{mod}\ 2).
\end{equation*} \qed

\begin{thm}\label{theo2}For any prime $p \geq 5$, $\alpha \geq 0$, and $n \geq 0$, we have
\begin{equation*}
b_2(p^{2\alpha+1}n+\frac{(24j+1)p^{2\alpha}-1}{24}) \equiv 0
\quad (\text{mod}\ 2),
\end{equation*}
where the integer $j$ satisfies $0 \leq j \leq p-1$ and $\left(\frac{24j+1}{p}\right)=-1$.
\end{thm}
\pf For any prime $p \geq 5$ and an integer $j$ with $0 \leq j \leq p-1$, according to Theorem \ref{f-p}
and Lemma \ref{b2-f}, if $j \not\equiv (3k^2+k)/2 \ (\text{mod}\ p)$ for $|k|
\leq (p-1)/2$, then we have
\begin{equation*}
b_{2}(p^{2\alpha}(pn+j)+\frac{p^{2\alpha}-1}{24})\equiv
0 \quad (\text{mod}\ 2).
\end{equation*}\qed

For example, we have the following congruences by setting $p=5$ in Theorem \ref{theo1} and Theorem \ref{theo2}.
\begin{exam} For $\alpha \geq 0$ and $n \geq 0$, we have
\begin{align*}
b_{2}(5^{2\alpha+2}n+\frac{29\cdot 5^{2\alpha+1}-1}{24}) \equiv 0 \quad (\text{mod}\ 2),\\
b_{2}(5^{2\alpha+2}n+\frac{53\cdot 5^{2\alpha+1}-1}{24}) \equiv 0 \quad (\text{mod}\ 2),\\
b_{2}(5^{2\alpha+2}n+\frac{77\cdot 5^{2\alpha+1}-1}{24}) \equiv 0 \quad (\text{mod}\ 2),\\
b_{2}(5^{2\alpha+2}n+\frac{101\cdot 5^{2\alpha+1}-1}{24}) \equiv 0 \quad (\text{mod}\ 2),\\
b_{2}(5^{2\alpha+1}n+\frac{73\cdot 5^{2\alpha}-1}{24}) \equiv 0 \quad (\text{mod}\ 2),\\
b_{2}(5^{2\alpha+1}n+\frac{97\cdot 5^{2\alpha}-1}{24}) \equiv 0 \quad (\text{mod}\ 2).
\end{align*}
\end{exam}

In fact, we can obtain the following more generalized congruences.

\begin{lem}\label{fff} For the primes $p_1, p_2,\ldots,p_r \geq 5$, $r \geq 0$, and $n \geq 0$, we have
\begin{equation*}
\sum_{n=0}^{\infty}b_2(\prod_{s=1}^{r}p_s^2n+\frac{\prod_{s=1}^{r}p_s^2-1}
{24})q^n \equiv f(-q) \quad (\text{mod}\ 2).
\end{equation*}
By convention, we set $\prod_{s=1}^{0}p_s^2=1$.
\end{lem}

\pf We prove it by induction on $r$. First, we get the initial case from Lemma \ref{b2-f}. Then suppose
\begin{equation*}
\sum_{n=0}^{\infty}b_2(\prod_{s=1}^{r}p_s^2n+\frac{\prod_{s=1}^{r}p_s^2-1}
{24})q^n \equiv f(-q) \quad (\text{mod}\ 2).
\end{equation*}
Based on Theorem \ref{f-p}, for a prime $p_{r+1} \geq 5$, we have
\begin{align*}
\sum_{n=0}^{\infty}b_2(\prod_{s=1}^{r}p_s^2(p_{r+1}^2n+\frac{p_{r+1}^2-1}{24})+\frac{\prod_{s=1}^{r}p_s^2-1}
{24})q^n&=\sum_{n=0}^{\infty}b_2(\prod_{s=1}^{r+1}p_s^2n+\frac{\prod_{s=1}^{r+1}p_s^2-1}
{24})q^n\\  &\equiv f(-q) \quad (\text{mod}\ 2).
\end{align*}
Therefore, the lemma holds for $r+1$.\qed

\begin{thm}\label{g-b2-1} For the primes $p_1, p_2,\ldots,p_r \geq 5$, $r \geq 1$, and $n \geq 0$, we have
\begin{equation*}
b_2(\prod_{s=1}^{r}p_s^2 n+\frac{(24i+p_{r})\prod_{s=1}^{r-1}p_s^2p_{r}-1}{24})
\equiv 0 \quad (\text{mod}\ 2),
\end{equation*}
where $i=1,2,\ldots, p_{r}-1$.
\end{thm}
\pf
Based on Theorem \ref{f-p} and Lemma \ref{fff}, we have
\begin{equation*}
\sum_{n=0}^{\infty}b_2(\prod_{s=1}^{r}p_s^2(p_{r+1}n+\frac{p_{r+1}^2-1}{24})+\frac{\prod_{s=1}^{r}p_s^2-1}
{24})q^n \equiv f(-q^{p_{r+1}}) \quad (\text{mod}\ 2).
\end{equation*}
Therefore, for $i=1,2,\ldots,p_{r+1}-1$, we have
\begin{equation*}
b_2(\prod_{s=1}^{r}p_s^2(p_{r+1}(p_{r+1}n+i)+\frac{p_{r+1}^2-1}{24})+\frac{\prod_{s=1}^{r}p_s^2-1}
{24}) \equiv 0 \quad (\text{mod}\ 2).
\end{equation*}
That means
\begin{equation*}
b_2(\prod_{s=1}^{r+1}p_s^2 n+\frac{(24i+p_{r+1})\prod_{s=1}^{r}p_s^2p_{r+1}-1}{24})
\equiv 0 \quad (\text{mod}\ 2).
\end{equation*} \qed

\begin{thm}\label{g-b2-2} For the primes $p_1, p_2,\ldots,p_r \geq 5$, $r \geq 1$, and $n \geq 0$, we have
\begin{equation*}
b_2(\prod_{s=1}^{r-1}p_s^2p_{r} n+\frac{(24j+1)\prod_{s=1}^{r-1}p_s^2-1}
{24}) \equiv 0 \quad (\text{mod}\ 2),
\end{equation*}
where the integer $j$ satisfies $0 \leq j \leq p_r-1$ and $\left(\frac{24j+1}{p_r}\right)=-1$.
\end{thm}
\pf According to Theorem \ref{f-p} and Lemma \ref{fff}, for the primes $p_1,p_2,\ldots,p_{r+1} \geq 5$, if the integer $j$ satisfies $0 \leq j
\leq p_{r+1}-1$ and $j \not\equiv (3k^2+k)/2 \ (\text{mod}\ p_{r+1})$ for $|k|
\leq (p_{r+1}-1)/2$, then we have
\begin{equation*}
b_2(\prod_{s=1}^{r}p_s^2(p_{r+1} n+j)+\frac{\prod_{s=1}^{r}p_s^2-1}
{24}) \equiv 0 \quad (\text{mod}\ 2).
\end{equation*}\qed

For example, we set $p_1=7$ and $p_2=5$ in Theorem \ref{g-b2-1} and Theorem \ref{g-b2-2}.
\begin{exam}We have
\begin{align*}
b_2(1225n+296)&\equiv 0\quad (\text{mod}\ 2),\\
b_2(1225n+541)&\equiv 0\quad (\text{mod}\ 2),\\
b_2(1225n+786)&\equiv 0\quad (\text{mod}\ 2),\\
b_2(1225n+1031)&\equiv 0\quad (\text{mod}\ 2),\\
b_2(245n+149)&\equiv 0\quad (\text{mod}\ 2),\\
b_2(245n+198)&\equiv 0\quad (\text{mod}\ 2).
\end{align*}
\end{exam}

\subsection{$4$-regular partitions and $13$-regular
partitions}

Based on the $p$-dissection identity of $\psi(q)$ in Theorem \ref{psi-p}, we have the following result.
\begin{lem}\label{b4-psi}
For any odd prime $p$, $\alpha \geq 0$, and $n \geq 0$, we have
\begin{equation*}
\sum_{n=0}^{\infty}b_4(p^{2\alpha}n+\frac{p^{2\alpha}-1}{8})q^n
\equiv \psi(q) \quad (\text{mod}\ 2).
\end{equation*}
\end{lem}
\pf We prove the lemma by induction on $\alpha$.

When $\alpha=0$, we know that
\begin{align*}
\sum_{n=0}^{\infty}b_4(n)q^n &= \frac{(q^4;q^4)_{\infty}}{(q;q)_{\infty}}\\
& \equiv (q;q)_{\infty}^3 \quad (\text{mod}\ 2)\\
&= \sum_{n=0}^{\infty}(-1)^n(2n+1)q^{\frac{n(n+1)}{2}}\qquad \text{by \eqref{Jacobi-3}}\\
&\equiv \sum_{n=0}^{\infty}q^{\frac{n(n+1)}{2}}\quad (\text{mod}\ 2)\\
&= \psi(q).
\end{align*}
Suppose the lemma holds for $\alpha$. Now we prove the case for
$\alpha+1$. According to Theorem \ref{psi-p}, we have
\begin{equation}\label{psi-pp}
\sum_{n=0}^{\infty}b_4(p^{2\alpha}(pn+\frac{p^2-1}{8})+\frac{p^{2\alpha}-1}{8})q^n
\equiv \psi(q^{p})\quad (\text{mod}\ 2).
\end{equation}
Moreover, we have
\begin{equation*}
\sum_{n=0}^{\infty}b_4(p^{2\alpha}(p^2n+\frac{p^2-1}{8})+\frac{p^{2\alpha}-1}{8})q^n
\equiv \psi(q)\quad (\text{mod}\ 2).
\end{equation*}
That is to say
\begin{equation*}
\sum_{n=0}^{\infty}b_4(p^{2\alpha+2}n+\frac{p^{2\alpha+2}-1}{8})q^n
\equiv \psi(q)\quad (\text{mod}\ 2).
\end{equation*}
We complete the proof.\qed

From the above lemma, we get the following two theorems.

\begin{thm}\label{b4-inf-1}For any odd prime $p$, $\alpha \geq 1$, and $n \geq 0$, we have
\begin{equation*}
b_4(p^{2\alpha}n+\frac{(8i+p)p^{2\alpha-1}-1}{8}) \equiv 0 \quad
(\text{mod}\ 2),
\end{equation*}
where $i=1,2,\ldots,p-1$.
\end{thm}

\begin{thm}\label{b444}For any odd prime $p$, $\alpha \geq 0$, and $n \geq 0$, we have
\begin{equation*}
b_4(p^{2\alpha+1}n+\frac{(8j+1)p^{2\alpha}-1}{8}) \equiv 0
\quad (\text{mod}\ 2),
\end{equation*}
where the integer $j$ satisfies $0 \leq j \leq p-1$ and $\left(\frac{8j+1}{p}\right)=-1$.
\end{thm}

Since the proofs of the above two theorems are similar to those of Theorem \ref{theo1} and
Theorem \ref{theo2} for $b_2(n)$, we omit the proofs.

Notice that we can get the congruences \eqref{b4-inf-23}-\eqref{b4-inf-25} given by Andrews et al. in \cite{Andrews-Hirschhorn-Sellers-2010}
by setting $p=3$ in Theorem \ref{b4-inf-1} and Theorem \ref{b444}, respectively. In addition, the above two theorems are
equivalent to the result given by Chen in \cite{Chen-2011}. While, similar to Theorem \ref{g-b2-1} and Theorem \ref{g-b2-2} for $b_2(n)$,
we can get more generalized congruences for $b_4(n)$.
\begin{thm} We have the following congrueces.
\begin{enumerate}
\item For the odd primes $p_1, p_2,\ldots,p_r$, $r \geq 0$, and $n \geq 0$, we have
\begin{equation*}
\sum_{n=0}^{\infty}b_4(\prod_{s=1}^{r}p_s^2n+\frac{\prod_{s=1}^{r}p_s^2-1}
{8})q^n \equiv \psi(q) \quad (\text{mod}\ 2).
\end{equation*}
By convention, we set $\prod_{s=1}^{0}p_s^2=1$.

\item For the odd primes $p_1, p_2,\ldots,p_r$, $r \geq 1$, and $n \geq 0$, we have
\begin{equation*}
b_4(\prod_{s=1}^{r}p_s^2 n+\frac{(8i+p_{r})\prod_{s=1}^{r-1}p_s^2p_{r}-1}{8})
\equiv 0 \quad (\text{mod}\ 2),
\end{equation*}
where $i=1,2,\ldots, p_{r}-1$.

\item For the odd primes $p_1, p_2,\ldots,p_r$, $r \geq 1$, and $n \geq 0$, we have
\begin{equation*}
b_4(\prod_{s=1}^{r-1}p_s^2p_{r} n+\frac{(8j+1)\prod_{s=1}^{r-1}p_s^2-1}
{8}) \equiv 0 \quad (\text{mod}\ 2),
\end{equation*}
where the integer $j$ satisfies $0 \leq j \leq p_r-1$ and $\left(\frac{8j+1}{p_r}\right)=-1$.
\end{enumerate}
\end{thm}
Here we omit the induction proof.

\begin{rem}
In the following discussion on the congruences for $b_{\ell}(n)$, we can also obtain
the above kind of congruences for different primes. While, in this paper, we focus on illustrating the type of
congruences with the same prime.
\end{rem}

In the next, by finding a congruence relation between the $4$-regular partition function and the $13$-regular partition function, we
obtain some new congruences for $b_{13}(n)$.

In \cite[therorem 2]{Calkin-Drake-James-2008}, the authors got
\begin{equation*}
\sum^{\infty}_{n=0}b_{13}(2n)q^{n}\equiv
(q^{2};q^{2})_{\infty}^{3}+q^{3}(q^{26};q^{26})_{\infty}^{3}\quad
(\text{mod}\ 2).
\end{equation*}
Then we find the fact:
\begin{equation}\label{13a}
\sum^{\infty}_{n=0}b_{13}(2(2n))q^{n}=\sum^{\infty}_{n=0}b_{13}(4n)q^{n}\equiv
(q;q)_{\infty}^{3}\equiv\sum^{\infty}_{n=0}b_{4}(n)q^{n}\quad
(\text{mod}\ 2).
\end{equation}
Therefore, based on Theorem \ref{b4-inf-1}, Theorem \ref{b444}, and the above congruence, we get the following theorem.
\begin{thm}For any odd prime $p$, $\alpha \geq 0$, and $n \geq 0$,
\begin{enumerate}
\item we have
\begin{equation*}
b_{13}(4 p^{2\alpha+2}n+\frac{(8i+p)p^{2\alpha+1}-1}{2}) \equiv 0 \quad
(\text{mod}\ 2),
\end{equation*}
where $i=1,2,\ldots,p-1$;

\item we have
\begin{equation*}
b_{13}(4 p^{2\alpha+1}n+\frac{(8j+1)p^{2\alpha}-1}{2}) \equiv 0
\quad (\text{mod}\ 2),
\end{equation*}
where the integer $j$ satisfies $0 \leq j \leq p-1$ and $\left(\frac{8j+1}{p}\right)=-1$.
\end{enumerate}
\end{thm}

Moreover, based on the congruence \eqref{133333} modulo $3$ for $b_{13}(n)$, we get the following results.

\begin{thm} For $\alpha \geq 1$ and  $n\geq 0$, we have
\begin{align*}
b_{13}(4\cdot3^{2\alpha+1}n+\frac{17\cdot 3^{2\alpha}-1}{2})&\equiv 0\quad (\text{mod}\ 6),\\
b_{13}(4\cdot3^{2\alpha}n+\frac{11\cdot 3^{2\alpha-1}-1}{2})&\equiv 0\quad (\text{mod}\ 6).
\end{align*}
\end{thm}
\pf For $\alpha\geq 1$, setting $\ell=2\alpha+1$ and $2\alpha$ in \eqref{133333}, respectively, we get two special cases
\begin{align*}
b_{13}(3^{2\alpha+1}(4n+2)+\frac{5\cdot 3^{2\alpha}-1}{2})
&= b_{13}(4\cdot3^{2\alpha+1}n+\frac{17\cdot 3^{2\alpha}-1}{2}) \equiv 0 \quad (\text{mod}\ 3),\\
b_{13}(3^{2\alpha}(4n+1)+\frac{5\cdot 3^{2\alpha-1}-1}{2})
&=b_{13}(4\cdot3^{2\alpha}n+\frac{11\cdot 3^{2\alpha-1}-1}{2}) \equiv 0 \quad (\text{mod}\ 3).
\end{align*}

Moreover, due to \eqref{b4-inf-23}, \eqref{b4-inf-24}, and \eqref{13a}, for $\alpha \geq 1$, we can get
\begin{align*}
b_{13}(4\cdot3^{2\alpha+1}n+\frac{17\cdot 3^{2\alpha}-1}{2})&\equiv 0\quad (\text{mod}\ 2),\\
b_{13}(4\cdot3^{2\alpha}n+\frac{11\cdot 3^{2\alpha-1}-1}{2})&\equiv 0\quad (\text{mod}\ 2).
\end{align*}
Therefore, we prove the theorem.
\qed


\subsection{$5$-regular partitions}

In \cite{Hirschhorn-Sellers-2010}, the authors got the following results.
\begin{align*}
&\sum_{n=0}^{\infty}b_{5}(n)q^{n}=\frac{(q^{8}; q^{8})_{\infty}(q^{20}; q^{20})^{2}_{\infty}}{(q^{2}; q^{2})_{\infty}^{2}(q^{40}; q^{40})_{\infty}}
+q\frac{(q^{4}; q^{4})^{3}_{\infty}(q^{10}; q^{10})_{\infty}(q^{40}; q^{40})_{\infty}}{(q^{2}; q^{2})_{\infty}^{3}(q^{8}; q^{8})_{\infty}(q^{20}; q^{20})_{\infty}},\\
&\sum_{n=0}^{\infty}b_{5}(2n)q^{n}\equiv (q^{2}; q^{2})_{\infty} \quad (\text{mod}\ 2),\\
&\sum_{n=0}^{\infty}b_{5}(2n+1)q^{n}\equiv (q^{10}; q^{10})_{\infty}\sum_{n=0}^{\infty}b_{5}(n)q^{n} \quad (\text{mod}\ 2).
\end{align*}
So we have
\begin{align*}
\sum_{n=0}^{\infty}b_{5}(2(2n)+1)q^{n}&=\sum_{n=0}^{\infty}b_{5}(4n+1)q^{n}\\
&\equiv (q^{5}; q^{5})_{\infty}\sum_{n=0}^{\infty}b_{5}(2n)q^{n} \quad (\text{mod}\ 2)\\
&\equiv (q^{5}; q^{5})_{\infty}(q^{2}; q^{2})_{\infty} \quad (\text{mod}\ 2).
\end{align*}

Based on the above fact, we get the following lemma.

\begin{lem}\label{b5f25}For any prime $p \geq 5$, $\left(\frac{-10}{p}\right)=-1$, $\alpha \geq 0$, and $n \geq 0$, we have
\begin{equation*}
\sum_{n=0}^{\infty}b_5(4 p^{2\alpha}n+\frac{7 p^{2\alpha}-1}{6})q^n\equiv f(-q^2)f(-q^5) \quad (\text{mod}\ 2).
\end{equation*}
\end{lem}
\pf We prove the lemma by induction on $\alpha$.
For $\alpha=0$, we know that
$$\sum_{n=0}^{\infty}b_{5}(4n+1)q^{n}\equiv f(-q^2)f(-q^5) \quad (\text{mod}\ 2).$$
Suppose the lemma holds for $\alpha$. Then for a prime $p\geq 5$, $\left(\frac{-10}{p}\right)=-1$, and $-(p-1)/2
\leq k,m \leq (p-1)/2$, we solve
$$2\cdot \frac{3k^2+k}{2}+5\cdot \frac{3m^2+m}{2} \equiv \frac{7p^2-7}{24}\quad (\text{mod}\ p).$$
That means
$$(12k+2)^2+10(6m+1)^2 \equiv 0  \quad (\text{mod}\ p).$$
Since $\left(\frac{-10}{p}\right)=-1$, there is only one solution $k=m=(\pm p-1)/6$ where $\pm$ depends on the condition
that $(\pm p-1)/6$ should be an integer. So there are no $k$ and $m$ for $-(p-1)/2
\leq k,m \leq (p-1)/2$ and $k,m \neq (\pm p-1)/6$ such that $(2\cdot (3k^2+k)/2+5\cdot (3m^2+m)/2)$ and
$(7p^2-7)/24$ are in the same residue class modulo $p$.

Therefore, using Theorem \ref{f-p}, we have
\begin{align*}
b_5(4 p^{2\alpha}(p^2n+\frac{7p^2-7}{24})+\frac{7 p^{2\alpha}-1}{6})&=b_5(4 p^{2\alpha+2}n+\frac{7 p^{2\alpha+2}-1}{6})\\
&\equiv f(-q^2)f(-q^5) \quad (\text{mod}\ 2).
\end{align*}
The lemma holds for $\alpha+1$.\qed

According to Lemma \ref{b5f25}, we get the following congruence for $b_5(n)$.

\begin{thm}\label{b5-2}For any prime $p \geq 5$, $\left(\frac{-10}{p}\right)=-1$, $\alpha \geq 1$, and $n \geq 0$, we have
\begin{equation*}
b_5(4 p^{2\alpha}n+\frac{(24i+7p)p^{2\alpha-1}-1}{6})\equiv 0 \quad (\text{mod}\ 2),
\end{equation*}
where $i=1,2,\ldots,p-1$.
\end{thm}
\pf Based on Lemma \ref{b5f25} and Theorem \ref{f-p}, we have
\begin{align*}
b_5(4 p^{2\alpha}(pn+\frac{7p^2-7}{24})+\frac{7 p^{2\alpha}-1}{6})&=b_5(4 p^{2\alpha+1}n+\frac{7 p^{2\alpha+2}-1}{6})\\
&\equiv f(-q^{2p})f(-q^{5p}) \quad (\text{mod}\ 2).
\end{align*}
Therefore, for $i=1,2,\ldots,p-1$, we have
\begin{equation*}
b_5(4 p^{2\alpha+1}(pn+i)+\frac{7 p^{2\alpha+2}-1}{6})\equiv 0 \quad (\text{mod}\ 2).
\end{equation*}
\qed

For example, by setting $p=17$, $i=3$, and $\alpha=1$, for $n \geq 0$, we have
$$b_5(1156n+541)\equiv 0 \quad (\text{mod}\ 2).$$

Note that we can get Theorem \ref{Hir-b5} given by Hirschhorn in \cite{Hirschhorn-Sellers-2010} by setting $\alpha=1$ in Theorem \ref{b5-2}.

\begin{lem}\label{kmj}For any prime $p\geq 5$ and $\left(\frac{-10}{p}\right)=-1$, given an integer $j$ with $0\leq j\leq p-1$,
there always exist integers $k$ and $m$ with $0\leq k,m\leq p-1$ such that $2\cdot (3k^{2}+k)/2+5 \cdot
(3m^{2}+m)/2\equiv j\ (\text{mod}\ p)$.
\end{lem}
\pf If \begin{equation*}
2\cdot \frac{3k^{2}+k}{2}+5 \cdot
\frac{3m^{2}+m}{2}\equiv j\ (\text{mod}\ p),
\end{equation*}
then \begin{equation*}
2(6k+1)^2+5(6m+1)^2\equiv 24j+7\quad
(\text{mod}\ p).
\end{equation*}
Let
$$x=6k+1, y=6m+1, a=24j+7.$$
Note that $\{0,1,\ldots,p-1\}$ is a complete residue system for $6k+1$ and $6m+1$ modulo $p$ for $0\leq k,m\leq p-1$.
Then we need to prove that there exist $x$ and $y$ such that
\begin{equation*}
2x^2+5y^2-a \equiv 0\quad
(\text{mod}\ p).
\end{equation*}

If $(2x^2+5y^2-a,p)=1$ for all $x$ and $y$, then according to the fact
$$\ell^{p-1}\equiv1\ (\text{mod}\ p), \text{ if } (\ell,p)=1,$$
we have
$$\sum_{x,y=0}^{p-1}(1-(2x^2+5y^2-a)^{p-1}) \equiv 0\quad(\text{mod}\ p).$$
Therefore, we only need to prove $\sum_{x,y=0}^{p-1}(1-(2x^2+5y^2-a)^{p-1}) \not\equiv 0\quad(\text{mod}\ p)$.

We have
\begin{align*}
&\sum_{x,y=0}^{p-1}(1-(2x^2+5y^2-a)^{p-1})\\
\equiv & -\sum_{x=0}^{p-1}\sum_{y=0}^{p-1}\sum_{m=0}^{p-1}{p-1 \choose m} (2x^2+5y^2)^m(-a)^{p-1-m}\quad(\text{mod}\ p)\\
= & -\sum_{m=0}^{p-1}{p-1 \choose m}(-a)^{p-1-m}\sum_{i=0}^m {m \choose i}2^i5^{m-i}\sum_{x=0}^{p-1}x^{2i}\sum_{y=0}^{p-1}y^{2m-2i}.
\end{align*}
According to the fact that for $k>0$
\begin{equation*}
\sum_{i=0}^{p-1}i^k\equiv \left\{
\begin{array}{lll}-1& (\text{mod}\ p),& \text{if } p-1\mid k,\\
0& (\text{mod}\ p),&\text{if } p-1 \nmid k,
\end{array}\right.
\end{equation*}
We have
\begin{equation*}
\sum_{x,y=0}^{p-1}(1-(3x^2+y^2-a)^{p-1})\equiv -{p-1 \choose \frac{p-1}{2}}2^{\frac{p-1}{2}} 5^{\frac{p-1}{2}}\not\equiv 0\quad(\text{mod}\ p).
\end{equation*}
Hence, we complete the proof.\qed

Based on Lemma \ref{kmj}, here we don't consider the similar congruences like those for $b_4(n)$ in Theorem \ref{b444}.

In the following, we use Ramanujan's identity \eqref{Ramanujan-qq} to get some more infinite family of congruences for $b_5(n)$.

For convenience, we set
\begin{equation*}
a(q)=\frac{(q^{10},q^{15};q^{25})_{\infty}}{(q^5,q^{20};q^{25})_{\infty}}\qquad \text{and} \qquad
b(q)=\frac{(q^5,q^{20};q^{25})_{\infty}}{(q^{10},q^{15};q^{25})_{\infty}}=\frac{1}{a(q)}.
\end{equation*}
Then, we can rewrite \eqref{Ramanujan-qq} as
\begin{equation}\label{Hir-qq}
(q;q)_{\infty}=(q^{25};q^{25})_{\infty}(a(q)-q-q^2b(q)).
\end{equation}

\begin{lem} \label{b5}For $\alpha \geq 0$ and $n\geq 0$, we have
\begin{equation*}
\sum^{\infty}_{n=0}b_{5}(4\cdot5^{2\alpha}n+\frac{7\cdot 5^{2\alpha}-1}{6})q^{n} \equiv f(-q^2)f(-q^5) \quad (\text{mod}\ 2).
\end{equation*}
\end{lem}
\pf We prove the congruence by induction on $\alpha$.

We know that the case for $\alpha=0$ holds since we have
$$\sum_{n=0}^{\infty}b_{5}(4n+1)q^{n}\equiv f(-q^2)f(-q^5) \quad (\text{mod}\ 2).$$

Suppose that the congruence holds for $\alpha$. By using \eqref{Hir-qq} on $f(-q^2)$, we have
\begin{align*}
\sum^{\infty}_{n=0}b_{5}(4\cdot5^{2\alpha}n+\frac{7\cdot 5^{2\alpha}-1}{6})q^{n}&\equiv (q^{2}; q^{2})_{\infty}(q^{5}; q^{5})_{\infty} \quad (\text{mod}\ 2)\\
&=(q^{5}; q^{5})_{\infty}(q^{50}; q^{50})_{\infty}(a(q^2)-q^{2}-q^{4}b(q^2)).
\end{align*}
Then we get
\begin{eqnarray*}
&&\sum^{\infty}_{n=0}b_{5}(4\cdot5^{2\alpha}(5n+2)+\frac{7\cdot 5^{2\alpha}-1}{6})q^{n}\\
&=&\sum^{\infty}_{n=0}b_{5}(4\cdot5^{2\alpha+1}n+\frac{11\cdot 5^{2\alpha+1}-1}{6})q^{n}\\
&\equiv& (q; q)_{\infty}(q^{10}; q^{10})_{\infty} \quad (\text{mod}\ 2)\\
&=& (q^{10}; q^{10})_{\infty}(q^{25}; q^{25})_{\infty}(a(q)-q-q^{2}b(q)) \qquad \text{by } \eqref{Hir-qq}.
\end{eqnarray*}
Therefore,
\begin{align*}
\sum^{\infty}_{n=0}b_{5}(4\cdot5^{2\alpha+1}(5n+1)+\frac{11\cdot 5^{2\alpha+1}-1}{6})q^{n}&=\sum^{\infty}_{n=0}b_{5}(4\cdot5^{2\alpha+2}n+\frac{7\cdot 5^{2\alpha+2}-1}{6})q^{n}\\
&\equiv (q^{2}; q^{2})_{\infty}(q^{5}; q^{5})_{\infty} \quad (\text{mod}\ 2).
\end{align*}
So the congruence holds for $\alpha+1$.\qed

\begin{thm} \label{b5m2}For $\alpha \geq 0$ and $n\geq 0$, we have
\begin{align*}
b_{5}(4\cdot5^{2\alpha+1}n+\frac{31\cdot 5^{2\alpha}-1}{6})&\equiv 0\quad (\text{mod}\ 2),\\
b_{5}(4\cdot5^{2\alpha+1}n+\frac{79\cdot 5^{2\alpha}-1}{6})&\equiv 0\quad (\text{mod}\ 2),\\
b_{5}(4\cdot5^{2\alpha+2}n+\frac{83\cdot 5^{2\alpha+1}-1}{6})&\equiv 0\quad (\text{mod}\ 2),\\
b_{5}(4\cdot5^{2\alpha+2}n+\frac{107\cdot 5^{2\alpha+1}-1}{6})&\equiv 0\quad (\text{mod}\ 2).
\end{align*}
\end{thm}
\pf According to the proof of Lemma \ref{b5}, we have
\begin{equation*}
\sum^{\infty}_{n=0}b_{5}(4\cdot5^{2\alpha}n+\frac{7\cdot 5^{2\alpha}-1}{6})q^{n}\equiv (q^{5}; q^{5})_{\infty}(q^{50}; q^{50})_{\infty}(a(q^2)-q^{2}-q^{4}b(q^2)) \quad (\text{mod}\ 2),
\end{equation*}
\begin{equation*}
\sum^{\infty}_{n=0}b_{5}(4\cdot5^{2\alpha+1}n+\frac{11\cdot 5^{2\alpha+1}-1}{6})q^{n} \equiv
(q^{10}; q^{10})_{\infty}(q^{25}; q^{25})_{\infty}(a(q)-q-q^{2}b(q))\quad (\text{mod}\ 2).
\end{equation*}
So we have
\begin{align*}
b_{5}(4\cdot5^{2\alpha}(5n+1)+\frac{7\cdot 5^{2\alpha}-1}{6}) & \equiv 0 \quad (\text{mod}\ 2),\\
b_{5}(4\cdot5^{2\alpha}(5n+3)+\frac{7\cdot 5^{2\alpha}-1}{6}) & \equiv 0 \quad (\text{mod}\ 2),\\
b_{5}(4\cdot5^{2\alpha+1}(5n+3)+\frac{11\cdot 5^{2\alpha+1}-1}{6}) &\equiv 0 \quad (\text{mod}\ 2),\\
b_{5}(4\cdot5^{2\alpha+1}(5n+4)+\frac{11\cdot 5^{2\alpha+1}-1}{6}) &\equiv 0 \quad (\text{mod}\ 2).
\end{align*}\qed

Note that we can get Theorem \ref{Calkin-b5} given by Calkin et al. in \cite{Calkin-Drake-James-2008} by setting $\alpha=0$ in the first two congruences of Theorem \ref{b5m2}.


\subsection{$8$-regular partitions and $16$-regular partitions}

\begin{lem}\label{b8-m}For any prime $p\equiv -1\ (\text{mod}\ 6)$, $\alpha\geq 0$, and $n \geq 0$, we have
\begin{equation*}
\sum_{n=0}^{\infty}
b_{8}(p^{2\alpha}n+\frac{7p^{2\alpha}-7}{24})q^{n} \equiv
\psi(q)f(-q^{4}) \quad(\text{mod}\ 2).
\end{equation*}
\end{lem}
\pf First,
\begin{align*}
\sum_{n=0}^{\infty} b_{8}(n)q^{n} &=
\frac{(q^8;q^8)_\infty}{(q;q)_\infty}\\
&\equiv \frac{(q; q)_\infty^{8}}{(q; q)_\infty} \quad(\text{mod}\ 2)\\
&= (q;q)_\infty^{7} \\
&\equiv (q; q)_\infty^3(q^{4}; q^{4})_\infty \quad (\text{mod}\ 2)\\
&\equiv \psi (q)f(-q^{4})\quad (\text{mod}\ 2).
\end{align*}
According to Theorem \ref{psi-p} and Theorem \ref{f-p}, for any
prime $p\equiv -1\ (\text{mod}\ 6)$, we discuss the
congruence properties modulo $p$ for the following form
\begin{equation*}
 \frac{k^{2}+k}{2}+4 \cdot \frac{3m^{2}+m}{2},
\end{equation*}
where $0\leq k\leq (p-1)/2$ and $-(p-1)/2\leq m \leq
(p-1)/2$.

When $k=(p-1)/2$ and $m=(-p-1)/6$, we have
\begin{equation*}
\frac{k^{2}+k}{2}+4 \cdot \frac{3m^{2}+m}{2} = \frac{7p^{2}-7}{24}.
\end{equation*}
In addition, if we have
\begin{equation}\label{k4m}
 \frac{k^{2}+k}{2}+4\cdot \frac{3m^{2}+m}{2}\equiv \frac{7p^{2}-7}{24}\quad
(\text{mod}\ p),
\end{equation}
then
$$3(2k+1)^{2}+(12m+2)^{2}\equiv 0\quad(\text{mod}\ p).$$
Since $\left(\frac{-3}{p}\right)=-1$ for $p\equiv -1\ (\text{mod}\ 6)$, we have the only one solution
$2k+1\equiv 12m+2\equiv 0\ (\text{mod}\ p)$ for \eqref{k4m}. That means
$k=(p-1)/2$ and $m=(-p-1)/6$. So there are no other $k$ and $m$ such that
$(k^{2}+k)/2+4\cdot (3m^{2}+m)/2$ and $(7p^{2}-7)/24$ are
in the same residue class modulo $p$.

Therefore, we get
\begin{equation}\label{88}
 \sum_{n=0}^{\infty} b_{8}(pn+\frac{7p^{2}-7}{24})q^{n}\equiv \psi(q^{p})f(-q^{4p}) \quad(\text{mod}\ 2).
\end{equation}
Then,
\begin{equation*}
 \sum_{n=0}^{\infty} b_{8}(p^2n+\frac{7p^{2}-7}{24})q^{n}\equiv \psi(q)f(-q^4) \quad(\text{mod}\ 2).
\end{equation*}

Following this rule, we can prove the lemma by induction on
$\alpha$. Here we omit the induction proof.\qed

\begin{thm}\label{thmmm}For any $p\equiv -1\ (\text{mod}\ 6)$, $\alpha\geq 1$, and $n \geq 0$,
 we have
\begin{align*}
 b_{8}(p^{2\alpha}n+\frac{(24i+7p)p^{2\alpha-1}-7}{24})
&\equiv 0  \quad(\text{mod}\ 2),
\end{align*}
where $i=1,\ldots,p-1$.
\end{thm}
\pf According to Lemma \ref{b8-m} and \eqref{88}, for $\alpha\geq 0$, we can get
\begin{equation*}
\sum_{n=0}^{\infty}
b_{8}(p^{2\alpha+1}n+\frac{7p^{2\alpha+2}-7}{24})q^{n} \equiv
\psi(q^p)f(-q^{4p}) \quad(\text{mod}\ 2).
\end{equation*}
Therefore, for $i=1, \ldots, p-1$, we have
\begin{align*}
b_{8}(p^{2\alpha+1}(pn+i)+\frac{7p^{2\alpha+2}-7}{24})= b_{8}(p^{2\alpha+2}n+\frac{(24i+7p)p^{2\alpha+1}-7}{24})\equiv 0 \quad(\text{mod}\ 2).
\end{align*}
\qed

Similar to Lemma \ref{kmj}, we know the fact: for any $p\equiv -1\ (\text{mod}\ 6)$ and the integer $j$ with $0\leq j\leq p-1$,
there always exist the integers $k$ and $m$ with $0\leq k,m\leq p-1$ such that $(k^{2}+k)/2+4 \cdot
(3m^{2}+m)/2 \equiv j\ (\text{mod}\ p)$. Therefore, we don't consider the similar congruences for $b_8(n)$ like those in Theorem \ref{b444}.

Given an example for Theorem \ref{thmmm}, we set $p=5$.
\begin{exam} For $\alpha\geq 1$ and $n \geq 0$, we have
\begin{align*}
b_{8}(5^{2\alpha}n+\frac{59\cdot5^{2\alpha-1}-7}{24})
\equiv 0 \quad(\text{mod}\ 2),\\
b_{8}(5^{2\alpha}n+\frac{83\cdot5^{2\alpha-1}-7}{24})
\equiv 0 \quad(\text{mod}\ 2),\\
b_{8}(5^{2\alpha}n+\frac{107\cdot5^{2\alpha-1}-7}{24})
\equiv 0 \quad(\text{mod}\ 2),\\
b_{8}(5^{2\alpha}n+\frac{131\cdot5^{2\alpha-1}-7}{24})
\equiv 0 \quad(\text{mod}\ 2).
\end{align*}
\end{exam}

Following the same process for $b_8(n)$, we discuss the congruences for $b_{16}(n)$.

We have
\begin{align*}
\sum_{n=0}^{\infty} b_{16}(n)q^{n}&\equiv \frac{(q; q)_\infty^{16}}{(q; q)_\infty}\quad (\text{mod}\ 2)\\
&= (q; q)_\infty^{15}\\
&\equiv (q; q)_\infty^{3}(q^{4}; q^{4})_\infty^{3} \quad (\text{mod}\ 2)\\
&\equiv \psi(q)\psi(q^{4})\quad (\text{mod}\ 2).
\end{align*}
For the prime $p \equiv -1\ (\text{mod}\ 4)$ and $0 \leq k,m \leq  p-1$, we disucss
$$\frac{k^2+k}{2}+4 \cdot \frac{m^2+m}{2} \equiv \frac{5p^2-5}{8}\quad (\text{mod}\ p).$$
For
$$(2k+1)^2+(4m+2)^2 \equiv 0\quad (\text{mod}\ p),$$
since $\left(\frac{-1}{p}\right)=-1$ if $p \equiv -1\ (\text{mod}\ 4)$, we have the only one solution $k=m=(p-1)/2$.

Therefore, similar to the results for $b_8(n)$, we get the following theorem.
\begin{thm} \label{t16}For any prime $p\equiv -1\ (\text{mod}\ 4)$, $\alpha\geq 0$, and $n \geq 0$, we have
\begin{align*}
\sum_{n=0}^{\infty}
b_{16}(p^{2\alpha}n+\frac{5p^{2\alpha}-5}{8})q^{n} &\equiv
\psi(q)\psi(q^4) \quad(\text{mod}\ 2),\\
b_{16}(p^{2\alpha+2}n+\frac{(8i+5p)p^{2\alpha+1}-5}{8})
&\equiv 0  \quad(\text{mod}\ 2),\quad i=1, 2,\ldots, p-1.
\end{align*}
\end{thm}
\pf Since the proof is similar to the case for $b_8(n)$, we only give a sketch of the proof.

The first congruence can be proved by induction on $\alpha$.

From the first congruence, we can get
$$\sum_{n=0}^{\infty}b_{16}(p^{2\alpha}(pn+\frac{5p^2-5}{8})+\frac{5p^{2\alpha}-5}{8})q^n\equiv \psi(q^p)\psi(q^{4p}) \quad(\text{mod}\ 2).$$
Therefore, we have
$$b_{16}(p^{2\alpha}(p(pn+i)+\frac{5p^2-5}{8})+\frac{5p^{2\alpha}-5}{8})\equiv 0 \quad(\text{mod}\ 2),\quad i=1, 2,\ldots, p-1.$$
\qed

Similar to Lemma \ref{kmj}, we can also prove the fact: for any prime $p\equiv -1\ (\text{mod}\ 4)$, given an integer $j$ with $0\leq j\leq p-1$,
there always exist integers $k$ and $m$ with $0\leq k,m\leq p-1$, such that $(k^{2}+k)/2+4 \cdot
(m^{2}+m)/2 \equiv j \ (\text{mod}\ p)$.

We give an example for $b_{16}(n)$ by setting $p=3$ in Theorem \ref{t16}.

\begin{exam} For $\alpha\geq 1$ and $n \geq 0$, we have
\begin{align*}
b_{16}(3^{2\alpha}n+\frac{23\cdot 3^{2\alpha-1}-5}{8})
&\equiv 0  \quad(\text{mod}\ 2),\\
b_{16}(3^{2\alpha}n+\frac{31\cdot 3^{2\alpha-1}-5}{8})
&\equiv 0  \quad(\text{mod}\ 2).
\end{align*}
\end{exam}

\section{More congruences for some $\ell-$regular partition functions}
Three famous congruences for $p(n)$ given by Ramanujan in \cite{Ramanujan-1919,Ramanujan-2000} are stated as
follows.
\begin{align*}
p(5n+4)&\equiv 0 \quad (\text{mod}\ 5),\\
p(7n+5)&\equiv 0 \quad (\text{mod}\ 7),\\
p(11n+6)&\equiv 0 \quad (\text{mod}\ 11).
\end{align*}
Moreover, there are another two beautiful congruences for $p(n)$, namely,
\begin{align*}
p(25n+24)&\equiv 0 \quad (\text{mod}\ 25),\\
p(49n+47)&\equiv 0 \quad (\text{mod}\ 49).
\end{align*}
Since
$$\sum_{n=0}^{\infty}b_{\ell}(n)q^n=\frac{(q^{\ell};q^{\ell})_{\infty}}{(q;q)_{\infty}}
=(q^{\ell};q^{\ell})_{\infty}\sum_{n=0}^{\infty}p(n)q^n,$$
we get the following lemma.
\begin{lem}\label{bp}For $k\geq1$, we have
\begin{align}
b_{5k}(5n+4)&\equiv 0 \quad (\text{mod}\ 5),\label{p(5k)-m}\\
b_{7k}(7n+5)&\equiv 0 \quad (\text{mod}\ 7),\nonumber\\
b_{11k}(11n+6)&\equiv 0 \quad (\text{mod}\ 11),\nonumber\\
b_{25k}(25n+24)&\equiv 0 \quad (\text{mod}\ 25),\nonumber\\
b_{49k}(49n+47)&\equiv 0 \quad (\text{mod}\ 49).\nonumber
\end{align}
\end{lem}

According to Lemma \ref{bp}, Theorem \ref{b5m2}, Theorem \ref{theom1.4}, and Theorem \ref{fp-t2}, we can get many more congruences for
some $\ell$-regular partition functions.
\begin{thm} For $\alpha \geq 0$ and $n\geq 0$ , we have
\begin{align*}
b_{5}(4\cdot 5^{2\alpha+2}n+\frac{31\cdot 5^{2\alpha}-1}{6})\equiv b_{5}(4\cdot 5^{2\alpha+2}n+\frac{79\cdot 5^{2\alpha}-1}{6})&\equiv 0 \quad (\text{mod}\ 10),\\
b_{5}(4\cdot 5^{2\alpha+3}n+\frac{83\cdot 5^{2\alpha+1}-1}{6})\equiv b_{5}(4\cdot 5^{2\alpha+3}n+\frac{107\cdot 5^{2\alpha+1}-1}{6})&\equiv 0 \quad (\text{mod}\ 10),\\
b_{7}(7\cdot3^{2\alpha+2}n+\frac{35\cdot 3^{2\alpha+1}-1}{4})\equiv b_{7}(7\cdot3^{2\alpha+3}n+\frac{77\cdot 3^{2\alpha+2}-1}{4})&\equiv 0 \quad (\text{mod}\ 21),\\
b_{25}(5\cdot3^{2\alpha+3}n+5\cdot 3^{2\alpha+2}-1)&\equiv 0 \quad (\text{mod}\ 15),\\
b_{25}(25\cdot3^{2\alpha+3}n+50\cdot 3^{2\alpha+2}-1)&\equiv 0 \quad (\text{mod}\ 75),\\
b_{49}(7\cdot3^{3\alpha+3}n+14\cdot 3^{3\alpha+2}-2)&\equiv 0 \quad (\text{mod}\ 21),\\
b_{49}(49\cdot3^{3\alpha+3}n+98\cdot 3^{3\alpha+2}-2)&\equiv 0 \quad (\text{mod}\ 147),\\
b_{10}(45n+39)&\equiv 0 \quad (\text{mod}\ 15),\\
b_{22}(297n+259)&\equiv 0 \quad (\text{mod}\ 33),\\
b_{28}(189n+117)&\equiv 0 \quad (\text{mod}\ 21).
\end{align*}
\end{thm}
\pf According to \eqref{p(5k)-m}, we have
$$b_{5}(5n+4) \equiv 0 \quad (\text{mod}\ 5).$$
So we have
\begin{align*}
b_5(5(4\cdot 5^{2\alpha+1}n+\frac{31\cdot 5^{2\alpha-1}-5}{6})+4)
=b_5(4\cdot 5^{2\alpha+2}n+\frac{31\cdot 5^{2\alpha}-1}{6})&\equiv 0 \quad (\text{mod}\ 5),\\
b_5(5(4\cdot 5^{2\alpha+1}n+\frac{79\cdot 5^{2\alpha-1}-5}{6})+4)
=b_5(4\cdot 5^{2\alpha+2}n+\frac{79\cdot 5^{2\alpha}-1}{6})&\equiv 0 \quad (\text{mod}\ 5),\\
b_5(5(4\cdot 5^{2\alpha+2}n+\frac{83\cdot 5^{2\alpha}-5}{6})+4)
=b_5(4\cdot 5^{2\alpha+3}n+\frac{83\cdot 5^{2\alpha+1}-1}{6})&\equiv 0 \quad (\text{mod}\ 5),\\
b_5(5(4\cdot 5^{2\alpha+2}n+\frac{107\cdot 5^{2\alpha}-5}{6})+4)
=b_5(4\cdot 5^{2\alpha+3}n+\frac{107\cdot 5^{2\alpha+1}-1}{6})&\equiv 0 \quad (\text{mod}\ 5).
\end{align*}
Due to the congruences in Theorem \ref{b5m2}, we have
\begin{align*}
b_5(4\cdot 5^{2\alpha+1}(5n)+\frac{31\cdot 5^{2\alpha}-1}{6})
=b_5(4\cdot 5^{2\alpha+2}n+\frac{31\cdot 5^{2\alpha}-1}{6}) &\equiv 0 \quad (\text{mod}\ 2),\\
b_{5}(4\cdot5^{2\alpha+1}(5n)+\frac{79\cdot 5^{2\alpha}-1}{6})
=b_5(4\cdot 5^{2\alpha+2}n+\frac{79\cdot 5^{2\alpha}-1}{6})&\equiv 0\quad (\text{mod}\ 2),\\
b_{5}(4\cdot5^{2\alpha+2}(5n)+\frac{83\cdot 5^{2\alpha+1}-1}{6})
=b_{5}(4\cdot5^{2\alpha+3}n+\frac{83\cdot 5^{2\alpha+1}-1}{6})&\equiv 0\quad (\text{mod}\ 2),\\
b_{5}(4\cdot5^{2\alpha+2}(5n)+\frac{107\cdot 5^{2\alpha+1}-1}{6})
=b_{5}(4\cdot5^{2\alpha+3}n+\frac{107\cdot 5^{2\alpha+1}-1}{6})&\equiv 0\quad (\text{mod}\ 2).
\end{align*}

Therefore, we get the congruences for the $5$-regular partition function in the theorem.

Similarly, based on Lemma \ref{bp}, we have
\begin{align*}
b_{7}(7(3^{2\alpha+2}n+\frac{5\cdot 3^{2\alpha+1}-3}{4})+5)
=b_{7}(7\cdot3^{2\alpha+2}n+\frac{35\cdot 3^{2\alpha+1}-1}{4}) &\equiv 0 \quad (\text{mod}\ 7),\\
b_{7}(7(3^{2\alpha+3}n+\frac{11\cdot 3^{2\alpha+2}-3}{4})+5)
=b_{7}(7\cdot3^{2\alpha+3}n+\frac{77\cdot 3^{2\alpha+2}-1}{4}) &\equiv 0 \quad (\text{mod}\ 7),\\
b_{25}(5(3^{2\alpha+3}n+ 3^{2\alpha+2}-1)+4)
=b_{25}(5\cdot3^{2\alpha+3}n+5\cdot3^{2\alpha+2}-1) &\equiv 0 \quad (\text{mod}\ 5),\\
b_{25}(25(3^{2\alpha+3}n+ 2\cdot3^{2\alpha+2}-1)+24)
=b_{25}(25\cdot3^{2\alpha+3}n+50\cdot3^{2\alpha+2}-1) &\equiv 0 \quad (\text{mod}\ 25),\\
b_{49}(7(3^{3\alpha+3}n+2\cdot 3^{3\alpha+2}-1)+5)
=b_{49}(7\cdot3^{3\alpha+3}n+14\cdot3^{3\alpha+2}-2) &\equiv 0 \quad (\text{mod}\ 7),\\
b_{49}(49(3^{3\alpha+3}n+2\cdot 3^{3\alpha+2}-1)+47)
=b_{49}(49\cdot3^{3\alpha+3}n+98\cdot 3^{3\alpha+2}-2) &\equiv 0 \quad (\text{mod}\ 49),\\
b_{10}(5(9n+7)+4)=b_{10}(45n+39)&\equiv 0 \quad (\text{mod}\ 5),\\
b_{22}(11(27n+23)+6)=b_{22}(297n+259)&\equiv 0 \quad (\text{mod}\ 11),\\
b_{28}(7(27n+16)+5)=b_{28}(189n+117)&\equiv 0 \quad (\text{mod}\ 7).
\end{align*}
Due to Theorem \ref{theom1.4} and Theorem \ref{fp-t2}, we get
\begin{align*}
b_{7}(3^{2\alpha+2}(7n+2)+\frac{11\cdot 3^{2\alpha+1}-1}{4})
=b_{7}(7\cdot3^{2\alpha+2}n+\frac{35\cdot 3^{2\alpha+1}-1}{4}) &\equiv 0 \quad (\text{mod}\ 3),\\
b_{7}(3^{2\alpha+3}(7n+6)+\frac{5\cdot 3^{2\alpha+2}-1}{4})
=b_{7}(7\cdot3^{2\alpha+3}n+\frac{77\cdot 3^{2\alpha+2}-1}{4}) &\equiv 0 \quad (\text{mod}\ 3),\\
b_{25}(3^{2\alpha+3}(5n+1)+2\cdot 3^{2\alpha+2}-1)
=b_{25}(5\cdot3^{2\alpha+3}n+5\cdot 3^{2\alpha+2}-1) &\equiv 0 \quad (\text{mod}\ 3),\\
b_{25}(3^{2\alpha+3}(25n+16)+2\cdot 3^{2\alpha+2}-1)
=b_{25}(25\cdot3^{2\alpha+3}n+50\cdot 3^{2\alpha+2}-1) &\equiv 0 \quad (\text{mod}\ 3),\\
b_{49}(3^{3\alpha+3}(7n+4)+2\cdot 3^{3\alpha+2}-2)
=b_{49}(7\cdot3^{3\alpha+3}n+14\cdot 3^{3\alpha+2}-2) &\equiv 0 \quad (\text{mod}\ 3),\\
b_{49}(3^{3\alpha+3}(49n+32)+2\cdot 3^{3\alpha+2}-2)
=b_{49}(49\cdot3^{3\alpha+3}n+98\cdot 3^{3\alpha+2}-2) &\equiv 0 \quad (\text{mod}\ 3),\\
b_{10}(9(5n+4)+3)=b_{10}(45n+39)&\equiv 0 \quad (\text{mod}\ 3),\\
b_{22}(27(11n+9)+16)=b_{22}(297n+259)&\equiv 0 \quad (\text{mod}\ 3),\\
b_{28}(27(7n+4)+9)=b_{28}(189n+117)&\equiv 0 \quad (\text{mod}\ 3).
\end{align*}
Therefore, we complete the proof.\qed

\section{Concluding remarks}

In \cite{Sellers-2003}, Sellers studied the $p$-regular partitions with distinct parts, and
found a parity result for this kind of partitions. Let $b_p'(n)$ denote the number of
the $p$-regular partitions with distinct parts of $n$.
\begin{thm}\label{Sellers-p}\cite[Theorem 2.1]{Sellers-2003}
Let $p$ be a prime greater than $3$ and let $r$ be an integer between $1$ and $p-1$,
inclusively, such that $24r+1$ is a quadratic nonresidue modulo $p$. Then, for all
nonnegative integers $n$, $b'(pn+r)\equiv 0\ (\text{mod}\ 2)$.
\end{thm}

Since
\begin{equation*}
\sum_{n=0}^{\infty}b_p'(n)q^n=\frac{(-q;q)_{\infty}}{(-q^p;q^p)_{\infty}}
\equiv \frac{(q;q)_{\infty}}{(q^p;q^p)_{\infty}}\quad (\text{mod}\ 2),
\end{equation*}
Due to Theorem \ref{f-p}, Theorem \ref{Sellers-p} can be easily obtained. Moreover, for the prime $p \geq 5$,
we can get
\begin{equation*}
\sum_{n=0}^{\infty}b_p'(pn+\frac{p^2-1}{24})q^n\equiv \frac{(q^p;q^p)_{\infty}}{(q;q)_{\infty}}\quad (\text{mod}\ 2).
\end{equation*}
Then we obtain a congruent relation between $b_p'(n)$ and $b_p(n)$.
\begin{equation*}
b_p'(pn+\frac{p^2-1}{24})\equiv b_p(n)\quad (\text{mod}\ 2).
\end{equation*}
Therefore, we can study arithmetic properties modulo $2$ of $b_p'(n)$ from those of $b_p(n)$.

In the forthcoming papers, using the similar technique in this paper, we will discuss arithmetic of
some other kinds of partitions, such as the broken $k$-diamond partitions, the $k$ dots bracelet partitions,
and $t$-core partitions.

\smallskip

\noindent {\bf Acknowledgements:} The authors would like to thank Professor Weidong Gao,
Dongchun Han, and Hanbin Zhang for many fruitful discussions. This work was supported by the
National Natural Science Foundation of China and the PCSIRT Project
of the Ministry of Education.

\end{document}